\documentclass{amsart}

\usepackage{amssymb,amscd}

\theoremstyle{plain}

\newtheorem*{thmA}{Theorem A}
\newtheorem*{thmB}{Theorem B}
\newtheorem{thm}{Theorem}[section]
\newtheorem{lem}[thm]{Lemma}

\theoremstyle{definition}

\newtheorem{dfn}[thm]{Definition}
\newtheorem*{conjecture}{Conjecture}
\newtheorem{question}{Question}

\theoremstyle{remark}
\newtheorem{remark}{Remark}

\DeclareMathOperator{\Hom}{Hom}
\DeclareMathOperator{\disc}{disc}
\DeclareMathOperator{\cont}{cont}
\DeclareMathOperator{\ctensor}{\widehat\otimes}

\newcommand{\N}{\mathbb{N}}
\newcommand{\Z}{\mathbb{Z}}
\newcommand{\F}{\mathbb{F}}
\newcommand{\Q}{\mathbb{Q}}
\newcommand{\UU}{\mathcal{U}}
\newcommand{\noreq}{\trianglelefteq}
\newcommand{\supp}{\operatorname{supp}}
\newcommand{\Hd}[1]{H_{\disc}^{#1}}
\newcommand{\Hc}[1]{H_{\cont}^{#1}}
\newcommand{\hd}[1]{H^{\disc}_{#1}}
\newcommand{\hc}[1]{H^{\cont}_{#1}}
\newcommand{\factor}[2]{{\raise0.7ex\hbox{$#1$} \!\mathord{\left/
 {\vphantom {#1 {#2}}}\right.\kern-\nulldelimiterspace}
\!\lower0.7ex\hbox{${#2}$}}}

\title[Discrete and Continuous Cohomology Groups of a Pro-$p$ Group]
{Comparison of the Discrete and Continuous Cohomology Groups of a Pro-$p$ Group}

\author[G.A. Fern\'andez-Alcober]{Gustavo A. Fern\'andez-Alcober}
\address{Matematika Saila, Zientzia eta Teknologia Fakultatea,
Euskal Herriko Unibertsitatea, 48080 Bilbao, Spain}
\email{gustavo.fernandez@ehu.es}

\author[I. Kazachkov]{Ilya V. Kazachkov}
\address{Department of Mathematics and Statistics, McGill University,
805 Sherbrooke St. West, Montreal, Quebec H3A 2K6, Canada}
\email{kazachkov@math.mcgill.ca}

\author[V. Remeslennikov]{Vladimir N. Remeslennikov}
\address{Institute of Mathematics (Russian Academy of Science),
13 Pevtsova St., Omsk 644099, Russia}
\email{remesl@iitam.omsk.net.ru}

\author[P. Symonds]{Peter Symonds}
\address{School of Mathematics, University of Manchester, P.O.\ Box 88,
Manchester M60 1QD, United Kingdom}
\email{peter.symonds@manchester.ac.uk}

\thanks{The first author is supported by the Spanish Ministry of Science
and Education, grant MTM2004-04665, partly with FEDER funds, and by
the University of the Basque Country, grant UPV05/99. The third author is supported by the Russian Foundation for Basic Research, grant 05-01-00057.}

\begin{document}

\begin{abstract}
To the Memory of L.~D.~Faddeev.
\end{abstract}

\maketitle

\section{Introduction}
\label{intro}

Let $G$ be a profinite group and let $M$ be a discrete $G$-module, i.e.\ $M$ has
the discrete topology and $G$ acts on $M$ continuously.
The continuous, or Galois, cohomology groups of $G$ with coefficients in $M$ can be
defined via continuous cochains from $G\times \overset{q}{\cdots} \times G$ to $M$.
General references for the continuous cohomology of profinite groups are the classical
monograph of Serre \cite{serre}, and some chapters of the recent books \cite{rz} and
\cite{wil} on profinite groups.

Of course, a profinite group $G$ can also be regarded as an abstract
group, and then the corresponding cohomology groups are obtained by
considering all cochains, and not only the continuous ones.
Topologists prefer to think of this as giving $G$ the discrete
topology, so we speak of discrete cohomology groups in this case.
Thus for a discrete $G$-module $M$ we can consider both the
continuous and the discrete cohomology groups. In order to avoid
confusion, we denote them by $\Hc{q}(G,M)$ and $\Hd{q}(G,M)$,
respectively. Since the continuous cochains form a subgroup of the
group of all cochains, it readily follows that there are natural
homomorphisms $\varphi^q:\Hc{q}(G,M) \to \Hd{q}(G,M)$, which we call
the {\em comparison maps\/} between the continuous and discrete
cohomology groups.

Our interest is focused on pro-$p$ groups.
If $P$ is a pro-$p$ group, then the trivial $P$-module $\F_p$ plays a prominent role
in the study of the cohomology of $P$, since $\F_p$ is the only discrete simple $p$-torsion
$P$-module.
As a consequence, in order to determine the continuous cohomological dimension of $P$,
one has to look only at the cohomology groups of $P$ with coefficients in $\F_p$.
Let us write, for short, $\Hc{q}(P)$ and $\Hd{q}(P)$ instead of $\Hc{q}(P,\F_p)$ and
$\Hd{q}(P,\F_p)$, respectively.

The following question was posed to one of the authors by G.~Mislin.

\begin{question}
If $P$ is a pro-$p$ group, when is the comparison map $\varphi^q:\Hc{q}(P) \to \Hd{q}(P)$
an isomorphism?
\end{question}

According to Mislin \cite{mislin}, the answer to this question is
affirmative for every $q$ if $P$ is poly-$\Z_p$ by finite. If we
look at small values of $q$, it is obvious that $\varphi^0$ is
always an isomorphism. On the other hand, since $\Hd{1}(P)$
coincides with the group of all homomorphisms from $P$ to $\F_p$ and
$\Hc{1}(P)$ is the subgroup of continuous homomorphisms, it is clear
that $\varphi^1$ is an isomorphism if and only if every
$f\in\Hom(P,\F_p)$ is continuous. This happens if and only if every
maximal subgroup of $P$ is open, which in turn is equivalent to $P$
being finitely generated. (See Exercise 6(d) in Section I.4.2 of
\cite{serre}, and Exercise 4.6.5 of \cite{wil}.) Now we ask: for
which finitely generated pro-$p$ groups is also $\varphi^2:\Hc{2}(P)
\to \Hd{2}(P)$ an isomorphism? As proved by Sury \cite{sury}, this
is the case if $P$ is either soluble $p$-adic analytic or Chevalley
$p$-adic analytic. In \cite{sury} B.~Sury attributes the following conjecture to G.~Prasad.

\begin{conjecture}[G.~Prasad]
For every $p$-adic analytic group $P$, the comparison map
$\varphi^2:\Hc{2}(P) \to \Hd{2}(P)$ is an isomorphism.
\end{conjecture}

We also address the corresponding problem for the homology of a
profinite group $G$. In this case, we have to consider profinite
$G$-modules. Recall that, as an abstract group, the homology of $G$
with coefficients in $M$ can be defined as the homology of the chain
complex which is obtained by tensoring $M$ with a projective
resolution $X\to \Z$ of $\Z$ over $\Z G$. (See Section III.1 of
\cite{brown}.) The homology of $G$ as a profinite group is defined
similarly, but this time we need to use a projective resolution
$X'\to \widehat\Z$ of $\widehat\Z$ in the category of profinite
modules over the completed group ring $\widehat\Z[[G]]$, and the
complete tensor product $\ctensor$ instead of the usual tensor
product $\otimes$. (See Section 6.3 of \cite{rz} in this case.)
Hence, again we have discrete homology groups $\hd{q}(G,M)$ and
continuous homology groups $\hc{q}(G,M)$. By the Comparison Theorem
(Theorem III.6.1 of \cite{mac}), there is a chain map $X\to X'$,
which gives rise to a chain map $X\otimes_{\Z G} M\to
X'\ctensor_{\widehat \Z[[G]]} M$. This way we obtain natural
homomorphisms $\varphi_q:\hd{q}(G,M)\to \hc{q}(G,M)$ between the
corresponding homology groups, and we are interested in determining
when $\varphi_q$ is an isomorphism. As before, we concentrate on the
homology groups of a pro-$p$ group $P$ with coefficients in the
trivial module $\F_p$, which we denote by $\hd{q}(P)$ and
$\hc{q}(P)$.

Our first main theorem is the following, which will be proved in Section \ref{thmA}
together with some other general results.

\begin{thmA}
Let $P$ be a finitely generated pro-$p$ group.
If $P$ is not finitely presented then $\varphi^2:\Hc{2}(P) \to \Hd{2}(P)$ is not
surjective.
Furthermore, if $P$ is finitely presented then the following conditions are equivalent:
\begin{enumerate}
\item
$\varphi^2:\Hc{2}(P) \to \Hd{2}(P)$ is an isomorphism.
\vspace{3pt}
\item
$\varphi_2:\hd{2}(P) \to \hc{2}(P)$ is an isomorphism.
\end{enumerate}
\end{thmA}

Note that Theorem A is not in conflict with the conjecture above, since $p$-adic analytic
groups are finitely presented \cite[Proposition 12.2.3]{wil}.
The following question arises.

\begin{question}
\label{question fp}
Does there exist a finitely presented pro-$p$ group $P$ for which
$\varphi^2:\Hc{2}(P) \to \Hd{2}(P)$ is not an isomorphism?
\end{question}

The second cohomology groups $\Hd{2}(P)$ and $\Hc{2}(P)$ classify the extensions of $P$ by
$\F_p$: arbitrary extensions in the first case, and extensions which are again a pro-$p$
group in the second.
(See Section IV.3 of \cite{brown} in the discrete case and Section 6.8 of \cite{rz} in the
continuous case.)
Since $\varphi^2$ is not surjective by Theorem A, for every finitely generated but not
finitely presented pro-$p$ group $P$, there exists an extension
\[
\begin{CD}
1 @>>> \F_p @>>> G @>>> P @>>> 1
\end{CD}
\]
of abstract groups such that $G$ can not be given the structure of a pro-$p$ group.
A different, and more difficult, question is to construct such extensions explicitly.
That is the purpose of the second main result in this paper.

\begin{thmB}
For every prime number $p$, and every power $q$ of $p$, there is a explicit
construction of an abstract group $G$ with a normal subgroup $Z$ of order $q$
such that:
\begin{enumerate}
\item
The quotient $P=\factor{G}{Z}$ can be endowed with a topology that makes this group
a pro-$p$ group.
Furthermore, $G$ can be chosen so that $P$ is finitely generated.
\item
It is not possible to define a topology on $G$ so that it becomes a pro-$p$ group.
\end{enumerate}
\end{thmB}

Unfortunately, as we will see later, the groups we construct do not
yield a finitely presented quotient $P$, and Question \ref{question
fp} remains open. However, Sections \ref{nfg} and \ref{fg} are interesting on their own right, as they provide a machinery for constructing (counter)examples in a non-standard way. In particular, we hope that some modification of our
construction can be used to provide examples with $P$ finitely
presented. On the other hand, even though $\varphi^2:\Hc{2}(P) \to
\Hd{2}(P)$ is not an isomorphism for the groups of Theorem B, it is
not clear whether $\varphi_2:\hd{2}(P) \to \hc{2}(P)$ is an
isomorphism or not.

Let us explain briefly the main idea behind these constructions.
Recall that a group $G$ is called {\em monolithic\/} if the intersection $T$
of all its non-trivial normal subgroups is non-trivial.
The subgroup $T$ is then called the {\em monolith\/} of $G$.
For example, a nilpotent group is monolithic if and only if $Z(G)$ is
monolithic.
In particular, a finite $p$-group is monolithic if and only if it has a cyclic
centre.
On the other hand, an infinite pro-$p$ group is never monolithic, since its
open normal subgroups have trivial intersection.
As a consequence, if $G$ is an infinite monolithic group then it is
impossible to define a topology on $G$ that makes it a pro-$p$ group.
Thus in order to prove Theorem B it suffices to construct monolithic groups
having a pro-$p$ group as a quotient.
In fact, we will obtain groups in which the centre is a cyclic $p$-group
and it is precisely the central quotient that will yield a pro-$p$ group.
More precisely, in Section \ref{nfg} we give a very general construction
of monolithic groups whose central quotient turns out to be the unrestricted
direct product of countably many copies of a certain finite $p$-group
(thus this quotient is an infinitely generated pro-$p$ group).
Important ingredients of this construction are central products, inverse limits
and ultraproducts.
Then, in Section \ref{fg}, we modify this construction with the help of
the wreath product in order to get finitely generated examples.

\section{General results and proof of Theorem A}
\label{thmA}

In this section, we prove Theorem A and some related results.
First of all, we collect some tools from the (co)homology of profinite groups
that will be needed.
In the following, let $G$ be a profinite group, and let $M$ be a topological
module on which $G$ acts continuously.
(Notice that, regardless of what kind of (co)homology we are considering, discrete
or continuous, the action of $G$ on all $G$-modules is assumed to be continuous.)

If $M$ is discrete, the continuous cohomology of $G$ with coefficients in $M$
can be obtained from the cohomology of its finite quotients
\cite[Corollary 6.5.6]{rz}:
\[
\Hc{q}(G,M) = \varinjlim \, H^q(\factor{G}{U},M^U),
\]
where $U$ runs over all open normal subgroups of $G$, and the direct limit
is taken over the inflation maps.
Here, $M^U$ is the subgroup of $U$-invariant elements of $M$.
As a consequence, if $P$ is a pro-$p$ group and $M$ a $p'$-module for $P$
(i.e.\ a $p'$-torsion $P$-module) then $\Hc{q}(P,M)=0$ for $q\ge 1$,
since the cohomology of a finite $p$-group with
coefficients in a $p'$-module is trivial for $q\ge 1$
\cite[Corollary III.10.2]{brown}.

In a similar fashion, if $M$ is profinite, the continuous homology
groups of $G$ with coefficients in $M$ are inverse limits
\cite[Corollary 6.5.8]{rz}:
\[
\hc{q}(G,M) = \varprojlim H_q(\factor{G}{U},M_U),
\]
where $M_U$ is the module of continuous coinvariants corresponding to $U$,
i.e.
\[
M_U = M/\overline{\langle gm-m \mid m\in M,\ g\in U \rangle},
\]
and the inverse limit is taken over the coinflation maps. It follows
that, for any pro-$p$ group $P$ and any $p'$-module $M$, the
homology groups $\hc{q}(P,M)$ are trivial for all $q\ge 1$.

As a consequence, we have the following description of the comparison maps
$\varphi^q$ and $\varphi_q$.

\begin{thm}
\label{limits}
For any profinite group $G$, the following hold:
\begin{enumerate}
\item
If $M$ is a discrete $G$-module, then the comparison map
$\varphi^q:\Hc{q}(G,M) \to \Hd{q}(G,M)$ is the direct limit of the inflation maps
$H^q(\factor{G}{U},M^U)\to \Hd{q}(G,M)$.
\item
If $M$ is a profinite $G$-module, then the comparison map
$\varphi_q:\hd{q}(G,M) \to \hc{q}(G,M)$ is the inverse limit of the
coinflation maps $\hd{q}(G,M) \to H_q(\factor{G}{U},M_U)$.
\end{enumerate}
\end{thm}

There are two dualising functors that concern us: $*_{\disc} =
\Hom_{\disc} (-,\factor{\Q}{\Z})$ and $*_{\cont} = \Hom_{\cont}
(-,\factor{\Q}{\Z})$. They take $G$-modules to $G$-modules, and
since we always deal with left $G$-modules, we are implicitly taking
the contragredient representation, i.e.\ changing the right action
into a left one by using the inverse map. Observe also that
$*_{\disc}$ and  $*_{\cont}$ agree on finite modules, in which case
we just write $*$. The second of these functors is Pontryagin
duality, which converts compact modules into discrete ones and
vice versa.
We know that $(M^{*_{\cont}})^{*_{\cont}} \cong M$ for
$M$ discrete torsion or profinite, although
$(M^{*_{\disc}})^{*_{\disc}} \cong M$ if and only if $M$ is finite.
(In this paper, when we speak of finite modules, we tacitly
understand that they are discrete or, what is the same in that case,
profinite.)

We have the following Duality Theorem
(see \cite[Proposition VI.7.1]{brown} and \cite[Proposition 6.3.6]{rz}).

\begin{thm}
\label{uct}
Let $G$ be a profinite group.
Then:
\begin{enumerate}
\item
$\Hd{q}(G,M^{*_{\disc}}) \cong \hd{q}(G,M)^{*_{\disc}}$, for every $G$-module $M$.
\item
$\Hc{q}(G,M^{*_{\cont}}) \cong \hc{q}(G,M)^{*_{\cont}}$, for every profinite $G$-module
$M$.
\end{enumerate}
\end{thm}

This result allows us to relate the homology comparison map
for a finite $G$-module with the cohomology comparison map
for its dual.

\begin{lem}
\label{ident}
Let $G$ be a profinite group and let $M$ be a finite $G$-module.
Then the map $\varphi^q:\Hc{q}(G,M^*) \to \Hd{q}(G,M^*)$ is equal to the composition
\begin{multline*}
\Hc{q}(G,M^*) \cong \hc{q}(G,M)^{*_{\cont}} \overset{\delta}{\longrightarrow}
\hc{q}(G,M)^{*_{\disc}}
\\
\overset{(\varphi_q)^*}{\longrightarrow} \hd{q}(G,M)^{*_{\disc}} \cong \Hd{q}(G,M^*),
\end{multline*}
where $\delta$ is the injection map and the isomorphisms come from the Duality
Theorem \ref{uct}.
\end{lem}

\begin{proof}
Every element of $\Hc{q}(G,M^*)$ is in the image of the inflation map from
$H^q(\factor{G}{U},(M^*)^U)$, for some open normal subgroup $U$ of $M$.
All the maps involved commute with inflation.
It is therefore sufficient to prove equality in the case when $G$ is finite.
But then all three maps $\varphi^q$, $\varphi_q$ and $\delta$ are the identity
and the result follows.
\end{proof}

\begin{thm}
\label{dual}
Let $G$ be a profinite group and let $M$ be a finite $G$-module.
Consider the comparison maps $\varphi^q:\Hc{q}(G,M^*) \to \Hd{q}(G,M^*)$ and
$\varphi_q:\hd{q}(G,M) \to \hc{q}(G,M)$.
Then:
\begin{enumerate}
\item
If $\varphi^q$ is surjective then $\varphi_q$ is injective, and if $\varphi_q$
is surjective then $\varphi^q$ is injective.
\item
Any two of the following conditions implies the third: $\varphi_q$ is an isomorphism,
$\varphi^q$ is an isomorphism, and $\hc{q}(G,M)$ is finite.
\end{enumerate}
\end{thm}

\begin{proof}
All the assertions follow easily from Lemma \ref{ident}.
If $\varphi^q$ is surjective, then $(\varphi_q)^*$ must be surjective, and so $\varphi_q$
is injective.
On the other hand, if $\varphi_q$ is surjective then $(\varphi_q)^*$ is injective, and
thus $\varphi^q$ is injective.
This shows that (i) holds.

Let us now prove (ii).
Since $\varphi^q$ is the composition of $\delta$, $\varphi_q$ and two isomorphisms, it is
clear that any two of these conditions implies the third: $\varphi_q$ is an isomorphism,
$\varphi^q$ is an isomorphism, and $\delta$ is an isomorphism.
Now, the injection $\delta:\hc{q}(G,M)^{*_{\cont}} \to \hc{q}(G,M)^{*_{\disc}}$ is an
isomorphism if and only if $\hc{q}(G,M)$ is finite.
\end{proof}

In what follows, we specialize to pro-$p$ groups.
Special attention is devoted to the trivial module $\F_p$, since it is the only discrete
simple $p$-torsion module for a pro-$p$ group.
As a consequence, if $P$ is a pro-$p$ group and $M$ a finite $P$-module of $p$-power order,
all composition factors of $M$ are isomorphic to $\F_p$.

The following result relates the cohomology of pro-$p$ groups in low dimension with group
presentations (see \cite[Sections I.4.2 and I.4.3]{serre}).

\begin{thm}
\label{dimensions}
Let $P$ be a pro-$p$ group.
Then:
\begin{enumerate}
\item
$\dim_{\F_p} \Hc{1}(P)$ is equal to the number of topological generators of $P$.
\item
If $P$ is finitely generated, then $\dim_{\F_p} \Hc{2}(P)$ is equal to the number of relations
in a presentation of $P$ as a pro-$p$ group with minimum number of generators.
\end{enumerate}
\end{thm}

Let us now examine the behaviour of the comparison maps in low
dimensions.

\begin{thm}
\label{h0}
Let $P$ be a pro-$p$ group.
Then:
\begin{enumerate}
\item
$\varphi^0$ is an isomorphism for every discrete $P$-module.
\item
$\varphi_0$ is an isomorphism if either the module is finite or $P$ is finitely generated.
\end{enumerate}
\end{thm}

\begin{proof}
As already mentioned in the introduction, (i) is obvious.
On the other hand, by \cite[Section III.1]{brown} and \cite[Lemma 6.3.3]{rz}, we have
\[
\hd{0}(P,M)=\factor{M}{N}
\quad
\text{and}
\quad
\hc{0}(P,M)=\factor{M}{\overline N},
\]
where $N=\langle gm-m \mid m\in M,\ g\in P \rangle$.
Thus $\varphi_0$ is an isomorphism precisely when $\overline N=N$.
This is clearly the case if the module $M$ is finite.
Suppose then that $P$ is finitely generated, by elements $g_1,\ldots,g_n$,
say.
Since $M$ is profinite in the case of homology, $M$ is compact.
Consequently the sub-abelian-group $T=\sum _{i=1}^n (g_i-1)M$
is also compact, and thus it is closed in $M$ and is a
$P$-submodule. Since $P$ acts trivially on the quotient
$\factor{M}{T}$, it follows that $T$ contains $\overline N$, and
thus $\overline N=N$ also in this case.
\end{proof}

\begin{thm}
\label{h1}
For finitely generated pro-$p$ groups and finite modules,
$\varphi^1$ and $\varphi_1$ are isomorphisms.
\end{thm}

\begin{proof}
Let $P$ be a finitely generated pro-$p$ group and let $M$ be a finite $P$-module.
By \cite[Lemma 6.8.6]{rz}, $\hc{1}(P)$ is isomorphic to $\factor{P}{\Phi(P)}$, where $\Phi(P)$
is the Frattini subgroup of $P$.
Thus $\hc{1}(P)$ is finite.
Let $W$ be an arbitrary composition factor of $M$.
Then either $W$ is isomorphic to $\F_p$ and $\hc{1}(P,W)$ is finite, or $W$ is a
$p'$-module and $\hc{1}(P,W)=0$.
It follows that $\hc{1}(P,M)$ is also finite.
According to Theorem \ref{dual} (ii), it suffices to prove the result for $\varphi^1$.

We know that $\Hc{1}(P,M)$ and $\Hd{1}(P,M)$ classify continuous and discrete sections
of $M \rtimes P \to P$, respectively.
(See \cite[Lemma 6.8.1]{rz} and \cite[Proposition IV.2.3]{brown}.)
Since all abstract homomorphisms between finitely generated \mbox{pro-$p$} groups
are continuous \cite[Corollary 1.21]{ddms}, we conclude that $\varphi^1$ is an
isomorphism.
\end{proof}

We only need the following lemma in order to proceed to the proof of Theorem A.
This result is adapted from Exercise 1 in \cite[Section I.2.6]{serre}.

\begin{lem}
\label{equiv}
Let $P$ be a pro-$p$ group and let $n$ be an integer.
Then the following conditions are equivalent:
\begin{enumerate}
\item
$\varphi_q:\hd{q}(P,M) \to \hc{q}(P,M)$ is an isomorphism for $q\le n$ and a surjection
for $q=n+1$, for every finite $P$-module $M$ of $p$-power order.
\item
$\varphi_q:\hd{q}(P,M) \to \hc{q}(P,M)$ is an injection for $q\le n$, for every finite
$P$-module $M$ of $p$-power order.
\item
$\varphi_q:\hd{q}(P) \to \hc{q}(P)$ is an injection for $q\le n$ and a surjection for
$q=n+1$.
\end{enumerate}
There is also a cohomology version, which is obtained by interchanging everywhere
discrete and continuous, and injection and surjection.
\end{lem}

\begin{proof}[Proof of Theorem A]
Since $P$ is finitely generated, we know from Theorems \ref{h0} and \ref{h1} that
$\varphi_0$ and $\varphi_1$ are isomorphisms for finite coefficients.
It follows from Lemma \ref{equiv} that $\varphi_2$ is surjective for coefficients
in $\F_p$.
Then $\varphi^2$ is injective by Theorem \ref{dual}.

Suppose first that $P$ is not finitely presented.
By Theorem \ref{dimensions}, the cohomology group $\Hc{2}(P)$ is infinite.
If $\varphi^2$ is surjective for $\F_p$, then it is an isomorphism, and
consequently $\hc{2}(P)$ is finite, by combining (i) and (ii)
of Theorem \ref{dual}.
Since $\F_p^*\cong \F_p$, it follows from Theorem \ref{uct} that
$\Hc{2}(P)\cong \hc{2}(P)^*$.
Hence $\Hc{2}(P)$ is also finite, which is a contradiction.

Assume finally that $P$ is finitely presented.
By Theorem \ref{dimensions}, $\Hc{2}(P)$ is finite.
Since $|\hc{2}(P)|=|\Hc{2}(P)|$, the equivalence of (i) and (ii) follows
from Theorem \ref{dual}.
\end{proof}

As a consequence, we see that extensions of the sort constructed in this paper
must always exist for pro-$p$ groups that are finitely generated but not
finitely presented; the problem is to construct them explicitly.
It is still possible that $\varphi_2$ is an isomorphism for these groups.
This also explains our interest in finding a finitely presented example.

\begin{remark}
Bousfield \cite{bous} has shown that if $F$ is a free pro-$p$ group on at least two
generators then $\hd{2}(F) \oplus \hd{3}(F)$ is uncountable.
On the other hand, it is known \cite[Theorem 7.7.4]{rz} that $\Hc{q}(F)=0$ for $q\ge 2$,
and by Theorem \ref{uct}, also $\hc{q}(F)=0$.
Hence, at least one of the maps $\varphi_2$ and $\varphi _3$ is not an isomorphism.
\end{remark}

\begin{remark}
There is an interesting analogy with what happens when a compact Lie group is
considered as a discrete group (see \cite{milnor}, \cite{fm}).
\end{remark}

At this point, it is not difficult to give a proof of Mislin's result
mentioned in the introduction, Theorem \ref{poly-Z_p} below.
For this purpose, we study the comparison maps for the group $\Z_p$.

\begin{thm}
\label{z}
Let $M$ be a finite $\Z_p$-module.
Then the comparison maps $\varphi^q:\Hc{q}(\Z_p,M) \to \Hd{q}(\Z_p,M)$ and
$\varphi_q:\hd{q}(\Z_p,M) \to \hc{q}(\Z_p,M)$ are isomorphisms for all $q$.
\end{thm}

\begin{proof}
We already know that the comparison maps are isomorphisms in degrees
$0$ and $1$, so we assume $q\ge 2$ in the remainder of the proof.
In that case, we are going to prove that $\Hc{q}(\Z_p,M)=\hd{q}(\Z_p,M)=0$.
It then follows from the Duality Theorem~\ref{uct} that also
$\hc{q}(\Z_p,M)=\Hd{q}(\Z_p,M)=0$, and $\varphi^q$ and $\varphi_q$ are
trivially isomorphisms.
Since any finite module is the direct sum of a $p$-module and a $p'$-module,
we may assume that either $M$ is a $p$-module or a $p'$-module.

Let us first deal with $\Hc{q}$.
If $M$ is a $p$-module then $\Hc{q}(\Z_p,M)=0$, since the cohomological
$p$-dimension of $\Z_p$ is $1$ \cite[Theorem 7.7.4]{rz}.
On the other hand, as mentioned at the beginning of this section,
$\Hc{q}(\Z_p,M)=0$ if $M$ is a $p'$-module.

For the case of $\hd{q}$, we use a result of Cartan (see \cite{cartan} or
Theorem V.6.4 in \cite{brown}), stating in particular that
$\hd{*}(G,R) \cong \bigwedge^*(G \otimes R)$ for any torsion-free abelian
group $G$ and any principal ideal domain $R$ on which $G$ acts trivially.
Here, $\hd{*}(G,R)$ is viewed as an $R$-algebra under the Pontryagin
product and $\bigwedge^*$ indicates the exterior algebra.
Since $q\ge 2$, it follows from this result that $\hd{q}(\Z_p,F)=0$ for every
field $F$ of prime order with trivial action of $\Z_p$.
As a consequence, $\hd{q}(\Z_p,M)=0$ if $M$ is either a $p$-module, or a
$p'$-module on which $\Z_p$ acts trivially, since in both cases all composition
factors of $M$ are cyclic of prime order with trivial action.

For a general $p'$-module, we can consider an open subgroup $U$ of $\Z_p$
acting trivially on $M$, since the action is continuous and $M$ is finite.
Then $U\cong \Z_p$ and $\hd{q}(U,M)=0$, as we have just shown.
Since $|\Z_p:U|$ is a $p$-power, this index is invertible in $M$ and we
deduce from \cite[Proposition III.10.1]{brown} that $\hd{q}(\Z_p,M)=0$.
This concludes the proof of the theorem.
\end{proof}

\begin{thm}
\label{poly-Z_p}
Let $P$ be a pro-$p$ group which is poly-$\Z_p$ by finite, and let $M$ be a
finite $P$-module.
Then $\varphi^q:\Hc{q}(P,M) \to \Hd{q}(P,M)$ and
$\varphi_q:\hd{q}(P,M) \to \hc{q}(P,M)$ are isomorphisms for all $q$.
\end{thm}

\begin{proof}
By hypothesis, there is a series $P=P_0\ge P_1\ge \cdots \ge P_r=1$ of
closed subgroups of $P$ such that $P_{i+1}\trianglelefteq P_i$ and
$\factor{P_i}{P_{i+1}}$ is either finite or isomorphic to $\Z_p$ for all $i$.
We argue by induction on $r$.
According to Theorem~\ref{z}, $\varphi^q$ is an isomorphism on each factor
$\factor{P_i}{P_{i+1}}$.
We can now assemble these isomorphisms by using the Lyndon-Hochschild-Serre
spectral sequence $H^r(\factor{G}{H},H^s(H)) \Rightarrow H^{r+s}(G)$ in both
cohomology theories and the Comparison Theorem for spectral sequences.
The proof for homology goes along the same lines.
\end{proof}

\section{Proof of Theorem B: examples with $P$ not finitely generated}
\label{nfg}

Let $E$ be a monolithic finite $p$-group. As mentioned in the
introduction, this means simply that $E$ has cyclic centre, so there
is certainly a vast choice for $E$. Just note that, for any
non-trivial power $q$ of $p$, it is possible to take $E$ so that its
centre has order $q$. For example, since the group of units
$\UU(\factor{\Z}{p^n\Z})$ has a cyclic subgroup of order $p^{n-1}$
if $p>2$ and $n\ge 2$ and of order $2^{n-2}$ if $p=2$ and $n\ge 3$,
it is easy to construct semidirect products of two cyclic $p$-groups
which have a cyclic centre of any desired order, and the group can be even of class
$2$ if we wish.

Let us write $E_n$ to denote the direct product of $n$ copies of $E$.
Inside the centre $Z_n$ of $E_n$ we consider the following subgroup:
\[
F_n
=
\{ (x^{(1)},\ldots,x^{(n)})\in Z_n \mid x^{(1)}\ldots x^{(n)}= 1 \}.
\]
Note that $F_n$ is generated by all tuples of the form
$(1,\ldots,1,x,1,\ldots,1,x^{-1},1,\ldots,1)$. Thus the quotient
$K_n=\factor{E_n}{F_n}$ can be viewed as the ``natural" central
product of $n$ copies of $E$, and it has the following property.

\begin{thm}
\label{central monolithic} The group $K_n$ is monolithic, with
centre $L_n=\factor{Z_n}{F_n}\cong Z(E)$.
\end{thm}

\begin{proof}
Since $Z(E)$ is cyclic, it suffices to prove the second assertion or, equivalently,
that if $x=(x^{(i)})\in E_n$ satisfies that $[x,E_n]\le F_n$ then $x^{(i)}\in Z(E)$
for all $i=1,\ldots,n$.
This follows easily arguing by way of contradiction.
\end{proof}

Consider now the projection map
\[
\begin{matrix}
\varphi_n & : & E_n & \longrightarrow & E_{n-1}
\cr
& & (x^{(1)},\ldots,x^{(n)}) & \longmapsto & (x^{(1)},\ldots,x^{(n-1)}).
\end{matrix}
\]
Note that $\varphi(Z_n)=Z_{n-1}$, so we get a surjective
homomorphism $\psi_n:\factor{E_n}{Z_n}\longrightarrow
\factor{E_{n-1}}{Z_{n-1}}$. Since there is a natural isomorphism
$\factor{K_n}{L_n}\cong \factor{E_n}{Z_n}$, we may suppose that the
homomorphisms $\psi_n$ are defined from $\factor{K_n}{L_n}$ onto
$\factor{K_{n-1}}{L_{n-1}}$. Thus we can consider the projective
system $\{\factor{K_n}{L_n},\psi_n\}$ over $\N$. We then define
\[
K = \big\{ (k_n)\in {\textstyle \prod_{n\in\N}} \, K_n \mid
(k_nL_n)\in \lim_{\longleftarrow}\, \factor{K_n}{L_n} \big \}.
\]
Note that $L=\prod_{n\in\N} \, L_n$ is a subgroup of $K$.

Before giving the last step of our construction, we digress for a while in order to
introduce the concepts of non-principal ultrafilters, ultraproducts and ultrapowers.

\begin{dfn}
A {\em filter\/} over a non-empty set $I$ is a non-empty family $\UU$ of subsets
of $I$ such that:
\begin{enumerate}
\item
The intersection of two elements of $\UU$ also lies in $\UU$.
\item
If $P\in\UU$ and $P\subseteq Q$ then also $Q\in\UU$.
\item
The empty set does not belong to $\UU$.
\end{enumerate}
The filter $\UU$ is called {\em principal\/} if it consists of all supersets of a
fixed subset of $I$, and it is called an {\em ultrafilter\/} if it is maximal in
the set of all filters over $I$ ordered by inclusion.
\end{dfn}

Let $\UU$ be a non-principal ultrafilter over $\N$ in the remainder.
Then $\UU$ enjoys the following two properties (see \cite{bar} for the proofs):
\begin{enumerate}
\item[(P1)]
If $\N=P_1\cup \cdots \cup P_k$ is a finite disjoint union then there exists an
index $i$ such that $P_i\in\UU$ and $P_j\not\in\UU$ for all $j\ne i$.
\item[(P2)]
All cofinite subsets of $\N$ belong to $\UU$.
\end{enumerate}
Given a family $\{H_n\}_{n\in\N}$ of groups we can form the {\em ultraproduct\/} of
this family by the non-principal ultrafilter $\UU$: this is the quotient of the
unrestricted direct product $\prod_{n\in\N} \, H_n$ by the subgroup $D$ consisting
of all tuples whose support does not lie in $\UU$.
If all groups $H_n$ are equal to the same group $H$, we speak of the {\em ultrapower\/}
of $H$ by $\UU$.
The following property is an easy consequence of (P1) and will play a fundamental role
in our discussion:
\begin{enumerate}
\item[(P3)]
If $H$ is a finite group then the ultrapower of $H$ by $\UU$ is isomorphic to $H$.
\end{enumerate}

Now, in order to complete the construction of the group $G$ that will fulfil the
conditions of Theorem B, we consider the subgroup
\[
D = \{ (k_n)\in L \mid \supp \, (k_n)\not\in \UU \}
\]
of $L$ and define $G=\factor{K}{D}$.

\begin{thm}
\label{non-fg}
Let $G$ be defined as above.
Then:
\begin{enumerate}
\item
$G$ is nilpotent and $Z(G)=\factor{L}{D}\cong Z(E)$ is a cyclic
$p$-group. Thus $G$ is monolithic and cannot be given the structure
of a pro-$p$ group.
\item
$P=\factor{G}{Z(G)}$ is a pro-$p$ group.
\end{enumerate}
\end{thm}

\begin{proof}
(i) If $E$ has nilpotency class $c$, it is clear from our
construction that $G$ is nilpotent of class at most $c$. Let us now
see that $Z(G)=\factor{L}{D}$. The inclusion $\supseteq$ is clear.
In order to see the reverse inclusion, let $(k_n)\in K$ be such that
$[(k_n),K]\le D$ and let us see that $(k_n)\in L$, i.e.\ that
$k_n\in L_n$ for all $n$. By way of contradiction, suppose that
$k_s\not\in L_s$ for some $s$ and assume that $s$ is minimal with
this property. Write $k_s=aF_s$ and $a=(a^{(1)},\ldots,a^{(s)})$.
Then $a^{(i)}\not\in Z(E)$ for some index $i$. We claim that $i=s$.
Otherwise write $k_{s-1}=bF_{s-1}$ with
$b=(b^{(1)},\ldots,b^{(s-1)})$. Since
$\psi_s(k_sL_s)=k_{s-1}L_{s-1}$, we have $\varphi_s(a)\equiv b
\pmod{Z_{s-1}}$ and consequently $b^{(i)}\equiv a^{(i)}
\pmod{Z(E)}$. It follows that $b^{(i)}\not\in Z(E)$ and
$k_{s-1}\not\in L_{s-1}$, which is a contradiction with the
minimality of $s$. Thus $i=s$ and $a^{(s)}\not\in Z(E)$. Then there
exists an element $y^{(s)}\in E$ such that $[a^{(s)},y^{(s)}]\ne 1$.
Now we define
\[
y_n
=
\begin{cases}
(1,\overset{n}{\ldots},1),
& \text{if $n<s$,}
\cr
(1,\overset{s-1}{\ldots},1,y^{(s)},1,\overset{n-s}{\ldots},1),
& \text{if $n\ge s$,}
\end{cases}
\]
and $k'_n=y_nF_n$ for all $n\in\N$.
Then the tuple $(k'_n)$ lies in $K$.

Let us see that the support of $[(k_n),(k_n')]$ is $[s,+\infty)\cap\N$.
Once we prove this, it follows from (P2) that this commutator is not
in $D$, contrary to our assumption that $[(k_n),K]\le D$.
Suppose then $n\ge s$ and let us check that $[k_n,k'_n]\ne 1$.
For this purpose, write $k_n=xF_n$ and $x=(x^{(1)},\ldots,x^{(n)})$.
Arguing as above, we have $x^{(s)}\equiv a^{(s)} \pmod{Z(E)}$.
Hence $[x^{(s)},y^{(s)}]=[a^{(s)},y^{(s)}]\ne 1$ and
\[
[x,y_n]
=
(1,\overset{s-1}{\ldots},1,[x^{(s)},y^{(s)}],1,\overset{n-s}{\ldots},1)
\not\in F_n.
\]
It follows that $[k_n,k'_n]\ne 1$, as desired.

Finally, note that $Z(G)=\factor{L}{D}$ is nothing but the
ultrapower of $Z(E)$ by the non-principal ultrafilter $\UU$. Now
(P3) yields that $Z(G)\cong Z(E)$ is cyclic.

(ii) It follows from (i) that $\factor{G}{Z(G)}\cong \factor{K}{L}$.
Now the result follows from the isomorphism $\factor{K}{L}\cong
\displaystyle \lim_{\longleftarrow}\, \factor{K_n}{L_n}$, which is
an immediate consequence of the First Isomorphism Theorem if we
consider the natural projection of $K$ onto $\displaystyle
\lim_{\longleftarrow}\, \factor{K_n}{L_n}$.
\end{proof}

\section{Proof of Theorem B: examples with $P$ finitely generated}
\label{fg}

In this section we combine the ideas from Section \ref{nfg} with the wreath
product construction in order to obtain the groups claimed in Theorem B
having a (topologically) finitely generated quotient which is a pro-$p$
group.
First we are going to combine central and wreath products as follows.
Let $H$ be a group and let $C$ be a finite group of order $n$.
Consider the regular wreath product $W=H\wr C$, let $B=H^n$ be the
corresponding base group and let $Z=Z(H)^n$.
We define the following subgroup of $Z$:
\[
F
=
\{ (x_1,\ldots,x_n)\in Z \mid x_1\ldots x_n=1 \}.
\]
Then $F\noreq W$ and we call the quotient $\factor{W}{F}$ the
{\em wreath-central product\/} of $H$ and $C$.
Our next theorem states that this construction behaves well when $H$
is a monolithic finite $p$-group.
We need a lemma for this.

\begin{lem}
\label{wreath}
Let $H$ be a non-abelian group and let $C$ be a non-trivial finite group.
Let $W=H\wr C$ be the regular wreath product of $H$ and $C$, and let
$B$ and $Z$ be as above.
Then $[B,w]\not\le Z$ for all $w\in W\setminus B$.
\end{lem}

\begin{proof}
Write $w=dx$ with $1\ne d\in C$ and $x=(x_c)_{c\in C}\in B$.
Choose an element $h\in H\setminus Z(H)$ and let $b\in B$ be the tuple
having the element $h$ at the position corresponding to the identity of
$C$ and $1$ elsewhere.
Then the value of the commutator $[b,w]=b^{-1}b^w=b^{-1}(b^d)^x$ at the
position corresponding to $d$ is $h^{x_d}\not\in Z(H)$.
Thus $[b,w]\not\in Z$.
\end{proof}

\begin{thm}
Let $H$ be a non-abelian finite $p$-group which is monolithic and
let $J$ be the wreath-central product of $H$ with another finite
$p$-group $C$.
Then $Z(J)\cong Z(H)$ and hence $J$ is again monolithic.
\end{thm}

\begin{proof}
Keep the notation before the lemma and let $K=\factor{B}{F}$, which
is the natural central product of $|C|$ copies of $H$. Arguing as in
the proof of Theorem \ref{central monolithic}, we get
$Z(K)=\factor{Z}{F}\cong Z(H)$. Thus it suffices to show that
$Z(J)=Z(K)$. Notice that the inclusion $\subseteq$ follows from the
previous lemma. On the other hand, since $[Z,W]\le F$ by the
definition of the action of the wreath product, we also have
$Z(K)\le Z(J)$ and we are done.
\end{proof}

Let now $E$ be a finite non-abelian monolithic $p$-group.
According to the last theorem, the wreath-central product $J_n$ of $E$ with a
cyclic group $\langle \alpha_n \rangle$ of order $p^n$ is again monolithic.
Let us write $E_n$ for the direct product of $p^n$ copies of $E$ (not $n$ copies
as in the previous section), $Z_n=Z(E_n)$ and
\[
F_n = \{ (x^{(1)},\ldots,x^{(p^n)})\in Z_n \mid x^{(1)}\ldots x^{(p^n)}=1 \}.
\]
We know from Theorem \ref{central monolithic} that the quotient
$K_n=\factor{E_n}{F_n}$ is monolithic, with centre
$L_n=\factor{Z_n}{F_n}$. Note that $J_n=K_n \rtimes \langle \alpha_n
\rangle$ is a semidirect product with kernel $K_n$ and complement
$\langle \alpha_n \rangle$, and that
$\factor{J_n}{L_n}=\factor{K_n}{L_n} \rtimes \langle \alpha_n
\rangle$.

Now we further assume that $E$ has class $2$ and that $|Z(E)|=q$ is a previously
chosen power of $p$.
(Recall the beginning of Section \ref{nfg}.)
We consider the map $\varphi_n:E_n\longrightarrow E_{n-1}$ defined by the
following rule: if $x=(x^{(i)})\in E_n$ then the $i$-th component of the
image $\varphi_n(x)$ is given by the product
\[
\prod_{r\equiv i\hspace{-6pt}\pmod{\!p^{n-1}}} \, x^{(r)},
\]
where the factors appear in increasing order of their indices. Note
that $\varphi_n$ is not a group homomorphism, for this we would need
$E$ to be abelian. However, since $x\equiv y \pmod{Z_n}$ implies
that $\varphi_n(x)\equiv \varphi_n(y) \pmod{Z_{n-1}}$, $\varphi_n$
induces a map $\psi_n:\factor{E_n}{Z_n}\longrightarrow
\factor{E_{n-1}}{Z_{n-1}}$. Since $\factor{E_n}{Z_n}$ is abelian
(this is why we need $E$ to be of class $2$), $\psi_n$ is a
homomorphism. As in the previous section, due to the existence of a
natural isomorphism $\factor{K_n}{L_n}\cong \factor{E_n}{Z_n}$, we
may suppose that $\psi_n$ is defined from $\factor{K_n}{L_n}$ onto
$\factor{K_{n-1}}{L_{n-1}}$.

\begin{thm}
By means of the rule $\psi_n(\alpha_n)=\alpha_{n-1}$, the
homomorphism $\psi_n$ extends to a homomorphism
$\factor{J_n}{L_n}\longrightarrow \factor{J_{n-1}}{L_{n-1}}$.
\end{thm}

\begin{proof}
Let $x\in E_n$ and let $y$ be the conjugate of $x$ by $\alpha_n$,
so that $y^{(i)}=x^{(i-1)}$, where $i-1$ is taken modulo $p^n$
(between $1$ and $p^n$).
If $a=\varphi_n(x)$ and $b=\varphi_n(y)$ then
\begin{align*}
b^{(i)}
&=
\prod_{r\equiv i\hspace{-6pt}\pmod{\!p^{n-1}}} \, y^{(r)}
=
\prod_{r\equiv i\hspace{-6pt}\pmod{\!p^{n-1}}} \, x^{(r-1)}
\\
&\equiv
\prod_{r\equiv i-1\hspace{-6pt}\pmod{\!p^{n-1}}} \, x^{(r)}
=
a^{(i-1)} \pmod{Z(E)},
\end{align*}
which means that $bZ_n$ is the conjugate of $aZ_n$ by $\alpha_{n-1}$.
This already implies the result.
\end{proof}

Thus we can consider the projective systems
$\{\factor{K_n}{L_n},\psi_n\}$ and $\{\factor{J_n}{L_n},\psi_n\}$
over $\N$. We then define
\[
K = \big\{ (k_n)\in {\textstyle \prod_{n\in\N}} \, K_n \mid
(k_nL_n)\in \lim_{\longleftarrow}\, \factor{K_n}{L_n} \big \}
\]
and
\[
J = \big\{ (j_n)\in {\textstyle \prod_{n\in\N}} \, J_n \mid
(j_nL_n)\in \lim_{\longleftarrow}\, \factor{J_n}{L_n} \big \}.
\]
If $L=\prod_{n\in\N} \, L_n$, we have $L\le K\le J$.
As in the previous section, we consider a non-principal ultrafilter $\UU$ over
$\N$ and define the subgroup $D$ of $L$ consisting of all tuples
whose support is not in $\UU$.
Let $G=\factor{J}{D}$, $H=\factor{K}{D}$ and $Z=\factor{L}{D}$.
Note that $Z\cong Z(E)$ is a central cyclic subgroup of $G$ of order $q$.

\begin{thm}
\label{H monolithic}
The group $H$ is monolithic with centre $Z$.
\end{thm}

\begin{proof}
As in the proof of Theorem \ref{non-fg}, we suppose that $[(k_n),K]\le D$
but that $k_n\not\in L_n$ for some $n$.
We claim that the set $S=\{ n\in\N \mid k_n\not\in L_n \}$ is of the form
$[s,+\infty)\cap \N$ for some $s$.
For this purpose, it suffices to see that if $n\in S$ and $m\ge n$ then
also $m\in S$.
This is an immediate consequence of the definition of $K$: \ since
$(k_n)\in K$, we have
$(\psi_{n+1}\circ \cdots \circ \psi_m)(k_mL_m)=k_nL_n\ne L_n$ and
consequently $k_m\not\in L_m$.

Write $k_n=x_nF_n$ for all $n\in\N$. Since $k_s\not\in L_s$, it
follows that the tuple $x_s\in E_s$ has an element $v\not\in Z(E)$
at some position $i(s)$. Let us choose $w\in E$ such that $[v,w]\ne
1$. Now we construct recursively a sequence $\{i(n)\}_{n\ge s}$ such
that $i(n)\equiv i(n-1)\pmod{p^{n-1}}$ and the element at position
$i(n)$ of $x_n$ does not lie in the centraliser $C_E(w)$. Indeed,
assume that $i=i(n-1)$ is already chosen and suppose by way of
contradiction that there is no $r\equiv i\pmod{p^{n-1}}$ with the
required property. Then $x_n^{(r)}\in C_E(w)$ for all $r\equiv
i\pmod{p^{n-1}}$ and it follows from the definition of $\varphi_n$
that also $x_{n-1}^{(i)}\in C_E(w)$, which is a contradiction. We
define
\[
y_n
=
\begin{cases}
(1,\overset{n}{\ldots},1),
& \text{if $n<s$,}
\cr
(1,\overset{i(n)-1}{\ldots},1,w,1,\overset{n-i(n)}{\ldots},1),
& \text{if $n\ge s$,}
\end{cases}
\]
and $k'_n=y_nF_n$ for all $n\in\N$.
Then the tuple $(k'_n)$ lies in $K$ and from this point onwards we can mimic
the proof of Theorem \ref{non-fg} to get the final contradiction that
$[(k_n),(k'_n)]\not\in D$.
\end{proof}

We are now ready to prove Theorem B in the finitely generated case.
If $P$ is a finitely generated (topological) group, let us write $d(P)$ to
denote the minimum number of (topological) generators of $P$.

\begin{thm}
\label{last}
The group $G$ is monolithic and $P=\factor{G}{Z}$ is a finitely
generated pro-$p$ group. More precisely, $d(P)=d(E)+1$.
\end{thm}

\begin{proof}
By considering the natural projection of $J$ onto $\displaystyle
\lim_{\longleftarrow} \, \factor{J_n}{L_n}$ it follows that
$\factor{J}{L}$ is isomorphic to this inverse limit. Thus
$P=\factor{G}{Z}\cong \factor{J}{L}$ is a pro-$p$ group. Also, since
$\factor{J_n}{L_n}\cong \factor{E}{Z(E)}\wr C_{p^n}$, we have
$d(\factor{J_n}{L_n})=d(E)+1$ for all $n$ and consequently also
$d(P)=d(E)+1$. (The last assertion follows basically from
Proposition 4.2.1 in \cite{wil}, check its proof.)

Let us see now that $\factor{G}{H}\cong \Z_p$.
Note that $\psi_n$ induces an epimorphism
$\overline\psi_n:\factor{J_n}{K_n}\longrightarrow
\factor{J_{n-1}}{K_{n-1}}$, hence we can consider the projective
system $\{\factor{J_n}{K_n},\overline\psi_n\}$. Since
$\factor{J_n}{K_n}=\langle \overline \alpha_n \rangle$ and $\psi_n$
maps $\alpha_n$ to $\alpha_{n-1}$, we have $\displaystyle
\lim_{\longleftarrow} \, \factor{J_n}{K_n} \cong \Z_p$. Since the
intersection of the kernels of the natural homomorphisms
$\pi_n:J\longrightarrow \factor{J_n}{K_n}$ is precisely $K$, it
follows that $\factor{G}{H}\cong \factor{J}{K}\cong \Z_p$, as
claimed.

Finally, let us prove that $G$ is monolithic. Since we already know
that $H$ is monolithic, it suffices to see that $N\cap H\ne 1$ for
every non-trivial normal subgroup $N$ of $G$. Obviously, we may
suppose that $N$ is not contained in $H$. Since $\factor{G}{H}\cong
\Z_p$ is abelian, we deduce that $[G,N]\le H\cap N$. Assume by way
of contradiction that $[G,N]=1$. If we write $N=\factor{M}{D}$ with
$D<M\not\le K$, we obtain that $[J,M]\le D$. Let $W_n$ be the wreath
product $E\wr \langle \alpha_n \rangle$, i.e. $W_n=E_n\rtimes
\langle \alpha_n \rangle$. Note that $J_n=\factor{W_n}{F_n}$. Then
to each subgroup of $J$ there corresponds a subgroup in $W_n$, by
first projecting onto $J_n$ and then taking its preimage in $W_n$.
Of course, the subgroup corresponding to $J$ is $W_n$. Let us call
$M_n$ and $T_n$ the subgroups corresponding to $M$ in $J_n$ and
$W_n$, respectively. Since $M$ is not contained in $K=\cap_{n\in\N}
\, \ker \pi_n$, we have $M_n\not\le K_n$ for some $n$.
Consequently, for that value of $n$, $T_n$ is not contained in the
base group $E_n$ of the wreath product $W_n$. It follows from Lemma
\ref{wreath} that $[W_n,T_n]\not\le Z_n$, while the condition
$[J,M]\le D$ implies that $[W_n,T_n]\le Z_n$. This contradiction
proves the result.
\end{proof}

In fact, the proof above shows that the subgroup $Z$ is the centre of $G$, so
that the quotient $\factor{G}{Z}$ which yields a pro-$p$ group is in fact the central
quotient of $G$.
Just note that from the condition $[G,N]=1$ we have deduced that $N\le H$.
Consequently $Z(G)\le Z(H)=Z$, as proved in Theorem \ref{H monolithic}.
The reverse inclusion is obvious.

Finally, observe that the pro-$p$ group $P=\factor{G}{Z}$ in the last
theorem is not finitely presented.
This is a consequence of the following result \cite[Corollary 12.5.10]{wil}:
if $P$ is a finitely presented soluble pro-$p$ group having a closed normal
subgroup $Q$ such that $\factor{P}{Q}\cong \Z_p$, then $Q$ is finitely generated.
In our case, $P$ is soluble and, as seen in the proof of Theorem \ref{last},
for the subgroup $Q=\factor{H}{Z}$ we have $\factor{P}{Q}\cong \factor{G}{H}\cong \Z_p$.
Since $Q\cong \factor{K}{L}\cong \varprojlim \, \factor{K_n}{L_n}$ is not finitely
generated, it follows that $P$ is not finitely presented.

\end{document}